\documentclass[12pt]{article}
\usepackage{amsmath}
\usepackage{amsfonts}
\usepackage{amssymb}
\usepackage{latexsym}
\usepackage{amscd}
\usepackage{color}

\input amssym.def
\newcommand{\qed}{\hfill \ensuremath{\Box}}

\newcommand{\norma}[1]{\| #1 \|}
\newcommand{\conj}[2]{\left \{ {#1} \, : \, {#2} \right \}}

\newcommand{\cvf}{\overset{\omega}{\rightarrow}}

\newcommand {\cvfe} {\overset{\omega^\ast}{\rightarrow}}

\newcommand {\R}{\mathbb{R}}

\newcommand {\N} {\mathbb{N}}

\newcommand{\sol}[1]{\text{sol}(#1)}

\newtheorem{df}{Definition}[section]

\newtheorem{prop}[df]{Proposition}
\newtheorem{teo}[df]{Theorem}
\newtheorem{corol}[df]{Corollary}
\newtheorem{lemma}[df]{Lemma}

\begin{document}

\title{\textsc{A note on the Banach lattice $c_0( \ell_2^n)$, its dual  and  its bidual}
\author{M. L. Louren\c co and V. C. C. Miranda \footnote{Supported by a CAPES PhD scholarship 88882.377946/2019-01} \\USP, S\~ao Paulo, Brazil}}
\date{}
\maketitle

\begin{abstract}
    The main purpose of this paper is to study some geometric and topological  properties on $ c_0$ sum of  the
finite dimensional Banach lattice $\ell_2^n$, its dual  and its  bidual. Among other results, we show that               
 the Banach lattices $c_0(\ell_2^n)$ has the strong Gelfand-Philips property, but does not have the positive Grothendieck property.
We also
prove that the closed unit ball of $l_{\infty}(\ell_2^n)$ is an almost limited set.


\end{abstract}

\section{Introduction}

Throughout this paper $X$ and $Y$ will denote Banach spaces, $E$ and $F$ will denote Banach lattices. We denote by $B_X$ the closed unit ball of $X$. In a Banach lattice, the additional lattice structure provides a large number of tools that are not available in more general Banach spaces.
This fact facilitates the study of geometric and topological properties of Banach lattices.  It is extremely important to add more examples of Banach lattices which satisfy or do not  some  geometric  or topological  properties. Our objective here is to study the Banach lattices  given by
 $\left ( \bigoplus_{n=1}^\infty \ell_2^n \right )_0, \, \left ( \bigoplus_{n=1}^\infty \ell_2^n \right )_1$ and $\left ( \bigoplus_{n=1}^\infty \ell_2^n \right )_\infty$ and
 describe which properties each of them satisfies.
 The importance of such Banach  spaces is due to the fact that the result presented by Stegall in {\cite{stegall}, where he showed that $\left ( \bigoplus_{n=1}^\infty \ell_2^n \right )_\infty$ does not have Dunford-Pettis property, but its predual, $\left ( \bigoplus_{n=1}^\infty \ell_2^n \right )_1,$  has it.

 We will start by recalling concepts of specific sets in Banach spaces and their consequences on their geometric  or topological properties.

A bounded set $A \subset X$ is \textit{Dunford-Pettis} (resp. \textit{limited}) if every weakly null sequence in $X'$ converges uniformly to zero on $A$ (resp. if every weak* null sequence in $X'$ converges uniformly to zero on $A$). Concerning these sets, we can consider a few properties in the class of Banach spaces.  A Banach space  $X$ has the
     \textit{DP property} (resp. DP* property) if every relatively weakly compact subset of $X$ is Dunford-Pettis (resp. if every relatively weakly compact subset of $X$ is limited). Or equivalently, $x_n'(x_n) \to 0$ for every $x_n \cvf 0$ in $X$ and $x_n' \cvf 0$ in $X',$ (resp. $x_n'(x_n) \to 0$ for every $x_n \cvf 0$ in $X$ and $x_n' \cvfe 0$ in $X'$). We say that $X$ has \textit{Gelfand-Phlips property} (in short \textit{GP property}), if every limited subset of $X$ is relatively compact.

Of course the DP* property implies the DP. On the other hand, $L_1[0,1]$ and $c_0$ are Banach spaces with the DP property without the DP*. Schur spaces have all three properties listed above.
 Separable and reflexive spaces are examples of Banach spaces with the GP property. For more information concerning those properties we refer \cite{alip, bourgdiest, car, drew1}.

In the class of Banach lattices, the lattice structure allows us to consider disjoint sequences. A sequence $(x_n) \subset E$ is \textit{disjoint} if $|x_n| \wedge |x_m| = 0$ for every $n \neq m$. A bounded subset $A \subset E$ is \textit{almost Dunford-Pettis} (resp. \textit{almost limited}) if every disjoint weakly null sequence in $E'$ converges uniformly to zero on $A$ (resp. if every disjoint weak* null sequence in $E'$ converges uniformly to zero on $A$). Next, we will give  some properties in Banach lattices that the definitions given for the sets above appear naturally. A Banach lattices $E$ has the \textit{weak DP property} (in short wDP) if every relatively weakly compact subset of $E$ is almost Dunford-Pettis. Or equivalently, if for all Banach space $Y$, every weakly compact operator $T: F \to Y$ is an \textit{almost Dunford-Pettis} operator,  that means, $T$ maps disjoint weakly null sequences of $F$ onto norm null sequences in $Y$.  We say that $E$ has the \textit{weak DP* property} (wDP*) if every relatively weakly compact subset of $E$ is almost limited, or equivalently, if $x_n'(x_n) \to 0$ for each weakly null sequence $(x_n) \subset F$ and each disjoint weak* null sequence $(x_n') \subset F'$. A Banach lattice $E$ is said to have the \textit{strong GP property} (sGP) if every almost limited subset of $E$ is relatively compact.

Of course the DP and the DP* properties imply, respectively, the wDP and the wDP*. In \cite{leung}, Leung gave the first example of a Banach lattice with the wDP and without the DP. In \cite{chen}, the authors showed that $L_1[0,1]$ has the wDP* property even though it does not have the DP*. Note that the sGP property is stronger than the GP. For instance, $L_1[0,1]$ does not have the GP property. We refer \cite{ardakani, bouras, chen} for more details concerning those properties.

Recall that a Banach space $X$ has the \textit{Grothendieck property} if every weak* null sequence in $X'$ is weakly null. For example, $\ell_\infty$ has the Grothendieck property. For a Banach lattice, we can consider the weak Grothendieck property and the positive Grothendieck property. From \cite{wnuk}, $E$ is said to have the \textit{positive Grothendieck property} if every positive weak* null sequence in $E'$ is weakly null. Following \cite{mach}, $E$ has the \textit{weak Grothendieck property} if every disjoint weak* null sequence in $E'$ is weakly null. Clearly, the Grothendieck property implies both the positive Grothendick and the weak Grothendieck. For instance, $\ell_1$ has the weak Grothendieck property, but it fails to have the positive Grothendieck, and $c$ is a Banach lattice with the positive Grothendieck property without the weak Grothendieck property.

\section{Results}

First, we are going to fix some notations.
Denote by $\ell_2^n$ the Banach lattice $\R^n$ with Euclidean norm and let
\begin{equation} \label{eq1}
    E = \left ( \bigoplus_{n=1}^\infty \ell_2^n \right )_0 = c_0(\ell_2^n), \quad E' = \left ( \bigoplus_{n=1}^\infty \ell_2^n \right )_1 = l_1(\ell_2^n) \quad \text{and} \quad E'' = \left ( \bigoplus_{n=1}^\infty \ell_2^n \right )_\infty=  l_{\infty}(\ell_2^n).
\end{equation}

We can consider in $E$ a natural structure of Banach lattice induced by its unconditional basis $(e_j ^i)_{i,j}$ where $e_j ^i = (0, \dots, 0, \overset{i}{\overbrace{e_j}}, 0, \dots )$ with $e_j = (0, \dots, 0, 1_{(j)}, 0, \dots)$. Thus $E'$ and $E''$ also are Banach lattices with their dual structures. Our goal in this section is to study the topological  properties of such Banach lattices.
In the following, $E$, $E'$ and $E''$ will be fixed as in (\ref{eq1}). It is known that
 $E'$ is a Schur space (see \cite{stegall}), and consequently, $E'$ has the DP property. However, its dual $E''$ does not have it. This was the first example of a Banach space with the DP property whose dual space does not have it.

Since $E'$ has the Schur property, it has the DP, the DP* and the sGP properties, then $E$ has the DP and (hence) wDP properties. In the next proposition, we will show that $E$ does not have the wDP* property.

\begin{prop} \label{wdp1}
The Banach lattice $E$ does not have the wDP* property.
\end{prop}

\noindent \textit{Proof: }
If $(e_n)$ is the Schauder basis in $c_0$ and if $T: c_0 \to E$ is the positive diagonal operator given by
$$ T(\alpha_j)_j = \begin{pmatrix}
\alpha_1 & 0  & 0  & \dots \\ 
 & \alpha_2 & 0 & \dots \\ 
 &  & \alpha_3 & \dots\\ 
 &  &  & 
\end{pmatrix}, $$
we have that $Te_n \cvf 0$ in $E$.
On the other hand, the sequence $e_{n,n}' = (0, \dots, 0, e_n, 0, \dots) \in E'$ is disjoint and weak* null with $e_{n,n}' (Te_n) = 1$ for every $n$. So $E$ does not have the wDP* property.
\qed

\

As  $E''$ does not have the DP property, it is natural to ask if   $E''$  has  the wDP property. And here we  will  show that $E''$ does not have it. To do this, we need the next two lemmas. 

\begin{lemma} \label{wdp2}
Let $F$ and $G$ be Banach lattices such that $F$ has the wDP property. If $T: F \to G$ is a surjective lattice isomorphism, then $G$ also has the wDP property.
\end{lemma}

\noindent \textit{Proof: }
Let $X$ be a Banach space and   $S: G \to X$ be a weakly compact operator. So $S\circ T: F \to X$ is a weakly compact operator. As  $F$  has the wDP property, then  $S\circ T$ is  an almost DP operator operator.   Now we show that  the operator  $S$    is an almost DP operator. Let $(y_n) \subset G$ be a disjoint weakly null sequence, so there exists $(x_n) \subset F$  a  disjoint   weakly null sequence  such that $Tx_n = y_n$ for all $n$. Therefore  $S(y_n) = S(T(x_n)$  and we have that  $S(y_n) \to 0$ in $X,$ so $S$ is an almost DP and the result follows
\qed

\

Let $F$ be a Banach lattice and $G \subset F$ a sublattice. We say thay $G$ is {\it  a complemented sublattice} of $F$ if there is  a bounded projection $P: F \to F.$ such that $P(F) = G.$

\begin{lemma} \label{wdp3}
Let $F$ a Banach lattice and let $G$ be a complemented sublattice of $F.$ 
If $F$ has the wDP property, then $G$ has wDP property.
\end{lemma}

\noindent \textit{Proof: }
Consider $X$ be a Banach space and   $T: G \to X$ be  a weakly compact operator. Let    $P: F \to F$  a bounded projection   such that $P(F) = G.$  So $T\circ P: F \to X$ is  also a weakly compact operator. Thus,  $T\circ P$ is an almost DP operator. If $(z_n)$ is a disjoint weakly null sequence in $G$, since $G$ is a sublattice of $F$, it follows that $(z_n)$ is a disjoint weakly null sequence in $F$. As  $T(z_n) = T(P(z_n))$  we get that $T(z_n)\to 0$ and the result is true.
\qed

\

Now we can prove that $E''$ does not have the wDP property.

\begin{prop} \label{wdp4}
The Banach lattice $E''$ does not have the wDP property.
\end{prop}

\noindent \textit{Proof: }
Consider the bounded linear operator 
$R : E' \to \ell_2$ given by
$$ R(x) = (x_{1,1} + x_{2,1} + \cdots, x_{2,2} + x_{3,2} + \cdots, \cdots ). $$
By Stegall {\cite{stegall},  we have that   $R': \ell_2 \to E''$ is an isomorphism  on $R'(\ell_2),$  and $R'(\ell_2)$ is a complemented subspace of $E",$ as Banach spaces.  It is easy to verify that $R'$ is a  lattice isomorphism on $R'(\ell_2)$  and  $R'(\ell_2)$ is a complemented sublattice of $E''.$
 Since $\ell_2$ does not have the wDP property, it follows from Lemma \ref{wdp2} that $R'(\ell_2)$ cannot have the wDP property. By Lemma \ref{wdp3} follows that $E''$ cannot have the wDP property.
\qed

\

 Since  every almost limited set is almost DP set, as a consequence of Proposition \ref{wdp4}, $E''$ cannot have the wDP* property. We observe  that $E$ has the GP property, because  $E$ is separable.

 Now, we will prove that $E$ has the sGP property. First we need the following Lemma.

\begin{lemma} \label{sgp1}
The Banach lattice $E$ is   Dedekind complete.
\end{lemma}

\noindent \textit{Proof: }
Let $A \subset E$ such that $a \leq x$ for every $a \in A$ and some $x \in E^+$. In particular,
$$ a = \begin{pmatrix}
a_{1,1} & a_{2,1} & a_{3,1} & \cdots \\ 
 & a_{2,2} & a_{3,2} & \cdots \\ 
 &  & a_{3,3} & \cdots \\ 
 &  &  & \cdots 
\end{pmatrix} \leq 
\begin{pmatrix}
x_{1,1} & x_{2,1} & x_{3,1} & \cdots \\ 
 & x_{2,2} & x_{3,2} & \cdots \\ 
 &  & x_{3,3} & \cdots \\ 
 &  &  & \cdots 
\end{pmatrix} = x $$
holds for every $a \in A$.
So, $a_{i,j} \leq x_{i,j}$ in $\R$ for every $i \in \N$ and $j = 1, \dots, i$.  As $\R$ is Dedekind complete, there are $z_{i,j} = \sup \conj{a_{i,j}}{a = (a_{k,l})_{l \leq k} \in A}$, for all  $i \in \N$  and $j = 1,\cdot, i.$ Now, let $z = (z_{i,j})_{j \leq i}$. Since $x_{i,j} \leq z_{i,j}$ holds for every $i \in \N$ and $j = 1, \dots, i$, it follows that $z \in E$. Now we prove that $z = \sup A$. In fact, if $y \in E$ is such that $a \leq y$ for every $a \in A$, so $a_{i,j} \leq y_{i,j}$ for every $i \in \N$ and $j = 1, \dots, i$. Thus $z_{i,j} \leq y_{i,j}$ for every $i \in \N$ and $j = 1, \dots, i,$, 
hence $z$ is the supremum of $A$ in $E$.
\qed

\

As a consequence of above lemma, we have  that $E$  has order continuous norm. We will conclude that $E$ has the sGP property showing that $E$ is a discrete Banach lattice. First, we recall that an element $x$ belongs
to a Banach lattice $F$ is \textit{discrete} if $x > 0$ and $|y| \leq x$ implies $y = tx$ for some real
number $t$. If every order interval $[0, y]$ in $F$ contains a discrete element, then $F$
is said to be a \textit{discrete Banach lattice}.

\begin{teo} \label{sgp2}
The Banach lattice $E$ has the sGP property.
\end{teo}

\noindent \textit{Proof: }
By Theorem 2.1 of \cite{ardakani}, it suffices to prove that $E$ is a discrete Banach lattice. Let
$[0,y]$ be an order interval in $E,$ take
$$ y = \begin{pmatrix}
y_{1,1} & y_{2,1} & y_{3,1} & \cdots \\
 & y_{2,2} & y_{3,2} & \cdots \\ 
 &  & y_{3,3} & \cdots \\ 
 &  &  & \cdots
\end{pmatrix} \quad \quad \text{and} \quad \quad x = \begin{pmatrix}
y_{1,1} & 0 & 0 & \cdots \\ 
 & 0 & 0 & \cdots \\ 
 &  & 0 &  \cdots \\ 
 &  &  &  \cdots
\end{pmatrix}. $$
So $x \in [0,y]$. If $|z| \leq x$ in $E$, it follows that
$$ z = \begin{pmatrix}
z_{1,1} & 0 & 0 & \cdots \\ 
 & 0 & 0 & \cdots \\ 
 &  & 0 & \cdots \\ 
 &  &  & \cdots
\end{pmatrix} $$
with $|z_{1,1}| \leq y_{1,1}$ in $\R.$ This implies that there exists a real number $t$ such that $z_{1,1} = t  y_{1,1}$. Thus $z = tx.$
The result follows.
\qed

\

We claim that $E''$ does not have the GP property. Indeed, consider the positive operator $S: \ell_\infty \to E''$ given by 
$$ S(\alpha_j)_j = \begin{pmatrix}
\alpha_1 & 0  & 0 & \dots \\ 
 & \alpha_2 & 0 & \dots \\ 
 &  & \alpha_3 & \dots\\ 
 &  &  & 
\end{pmatrix}. $$
As $(e_n)_n \subset \ell_\infty$ is a weakly null limited sequence, this implies that $(Se_n)$ is a weakly null limited sequence in $E''$ such that 
$\norma{Se_n}_\infty = 1$ for all $n$. That means  $E''$ does not have the Gelfand-Phillips property.

Next, we  study the Grothendieck type properties in $E$ and $E'$.

\begin{prop} \label{grothen1}
The Banach lattice $E$ does not have the weak Grothendieck property and  does not have the positive Grotendieck property
\end{prop}

\noindent \textit{Proof: }
Let $(e_{n,n}') \subset E'$ as given in the proof of Proposition \ref{wdp1}. This sequence is positive, disjoint and weak* null in $E',$ however, $(e_{n,n}')$ is not weakly null. Indeed, if  $x'' \in E''$ given by
$$x'' = \begin{pmatrix}
1 & 0 & 0 & 0 & \cdots \\ 
 & 1 & 0 & 0 & \cdots \\ 
 &  & 1& 0 & \cdots \\ 
 &  &  &0 & \cdots 
\end{pmatrix}.$$
then $x''(e_{n,n}') = 1$ for all $n$.
\qed

\begin{prop} \label{grothen2}
The Banach lattice $E'$ has the weak Grothendieck property, however it does not have the positive Grothendieck property.
\end{prop}

\noindent \textit{Proof: }
Let $e \in E'',$ given by

$$ e = \begin{pmatrix}
1 & 1 & 1 & \cdots \\ 
 & 1 & 1 & \cdots \\ 
 &  & 1 & \cdots \\ 
 &  &  & \cdots 
\end{pmatrix} $$
 It is easy to see that $e$ is an order unit  of  $E''$, that means  $B_{E''} = [-e,e]$.
Let $(x_n'') \subset E''$ be a disjoint weak* null sequence. In particular, $(x_n'')$ is bounded, and so  there exists $M > 0$ such that $x_n'' \in  [-M e, M e]$ for every $n \in \N$. Consequently, $(x_n'')$ is a disjoint order bounded sequence in $E'',$  hence $x_n'' \cvf 0$ in $E''$ (pg. 192, \cite{alip}). Therefore $E'$ has the weak Grothendieck property.

Consider  the diagonal operator $T: \ell_1 \to E',$ given by
$$ T(\alpha_j)_j = \begin{pmatrix}
\alpha_1 & 0  & 0  & \dots \\ 
 & \alpha_2 & 0 & \dots \\ 
 &  & \alpha_3 & \dots\\ 
 &  &  & 
\end{pmatrix}, $$
It is easy to see  thar $T$ is a lattice isometry and for each $\alpha = (\alpha_j)_j \in \ell_1$, we have $\Vert T (\alpha)\Vert =\Vert \alpha \Vert, $  hence  $T$ is a lattice embedding.  By the Proposition 2.3.11 of \cite{meyer} yields that $\ell_1$ is isomorphic to a positively complemented sublattice in $E'$, then $E'$ does not have the positive Grothendieck property.
\qed

\begin{corol}
\label{grothen3}
The norm in $E''$ is not order continuous.
\end{corol}
\noindent \textit{Proof: }
By Proposition 4.9 from \cite{mach}, if $F$ is a Banach lattice which has the weak Grothendieck property and $F'$ has order continuous norm, then $F$ also has the positive Grothendieck property. So, by 
Proposition \ref{grothen2}, it follows that $E''$ does not have order continuous norm.
\qed

\

Now we want  to give a version   of Phillip's Lemma for $E'',$ for this we use Dixmeir's Theorem, that is,
 if $X$ is a Banach space, then $X'$ is complemented in $X''$ (see \cite{dix}). Despite being a known result, we decided to state it in the next lemma and present a proof for the specific case of the Banach lattice $E'''.$

\begin{lemma} \label{grothen4}
 Consider  $E^\perp = \conj{f \in E'''}{f(x) = 0, \, \forall x \in E}.$  Then  $E''' = E' \oplus E^\perp$ and 
 $E^\perp$ is an ideal in $E'''$.

\end{lemma}

\noindent \textit{Proof: }
  Let $f \in E'''$ and put $a_{i,j} = f(e_{i,j})$ for all $i \in \N$ and $j = 1, \dots, i$. We claim that
$$ a' = \begin{pmatrix}
a_{1,1} & a_{2,1} & a_{3,1} & \cdots \\ 
 & a_{2,2} & a_{3,2} & \cdots \\ 
 &  & a_{3,3} & \cdots \\ 
 &  &  & \cdots
\end{pmatrix} \in E'. $$
Indeed, since 
$$ \sum_{i=1}^\infty \norma{a_i}_2 \leq \sum_{i=1}^\infty \sum_{j=1}^i |a_{i,j}| = \sum_{i=1}^\infty \sum_{j=1}^i f(\epsilon_{i,j} e_{i,j}) = f(\sum_{i=1}^\infty \sum_{j=1}^i \epsilon_{i,j} e_{i,j}) \leq \norma{f}.  $$
On the other hand, if $x = (x_{i,j})_{1 \leq j \leq i} \in E$, then
$$ f(x) - a'(x) = \lim_{n \to \infty} \left [\sum_{i=1}^n \sum_{j=1}^i f_{i,j}(x_{i,j} e_{i,j}) - \sum_{i = 1}^n \sum_{j=1}^i a_{i,j}(x_{i,j})\right ] = 0. $$
As $E' \cap E^\perp = \{0\}$, we get $E''' = E' \oplus E^\perp$.

  By Lemma 2.5  $E$ has order continuous norm, as a consequence  it is  an  ideal in $E''.$    We claim that $E^\perp$ is an ideal in $E'''.$
Indeed, let $x''', y''' \in E'''$ with $|y'''| \leq |x'''|$ and $x''' \in E^\perp,$
If $x \in E^+$ and $|y''| \leq x$ in $E''$, then $|y''| \in E$, this implies that
$|x'''|(x) = \sup \conj{|x'''(y)|}{|y| \leq x} = 0.$  Finally, if $x \in E$,
$$ |y'''(x)| \leq |y'''| (|x|) \leq |x'''|(|x|) = 0. $$
Therefore $y''' \in E^\perp$.
\qed

\

Now we give  a version of Phillip's Lemma for $E'',$ the proof   follows  the same idea   of Theorem 4.67  in  \cite{alip}. We remark that 
 in \cite{carlo} the authors has showed that   $B_E$ is a limited set in $E^{''},$ but they used  another technique in  another context.

\begin{prop} \label{grothen6}
Every weak* null sequence in $E''$ converges uniformly to zero on $B_E$. Consequently, $B_E$ is a limited set in $E''$.
\end{prop}
\noindent \textit{Proof: }
Let $(f_n) \subset E'''$
be a weak* null sequence. By  Lemma \ref{grothen4} we write $f_n = x_n + g_n$ with $(x_n) \subset E'$ and $(g_n) \subset E^\perp$. 
As $E'$ has order continuous norm then  $E'$ is an ideal in $E'''$. On the other hand, since $E^\perp$ also is an ideal in $E'''$, by Theorem 1.41 in \cite{alip}, we have that $E'$ is a projection band in $E'''$, which yields that $x_n \cvfe 0$ in $E'''$ (see Theorem 4.46 in \cite{alip}). Then $x_n \cvf 0$ in $E'$, and since $E'$ has the Schur property, $x_n \to 0$ in $E'$. As a consequence,
$ \norma{f_n}_{B_E} = \sup_{x \in B_E} |x_n(x)| \leq \norma{x_n} \to 0. $
\qed

\

\begin{prop} \label{grothen7}
 The  Banach lattice $E''$ has the  weak and  the positive Grothendieck properties. 
\end{prop}
\noindent \textit{Proof: }  Let  positive weak* null sequence $(x_n') \subset E'''.$  Since  $e$ is an order unit  of  $E''$, that means  $B_{E''} = [-e,e],$ we get that
$$\Vert  x _n' \Vert = \sup_{ x \in  [-e,e]} \vert x_n(x)' \vert = x_n'(e)   \to 0.$$
So  $x_n' \cvf 0,$ and $E''$ has the positive Grothendieck property.

Now,  as   every disjoint weak* null sequence $(f_n) \subset E'''$ implies $\vert f_n\vert  \cvfe  0,$ in $E''$  and $E''$  has the positive Grothendieck property the result follows.
\qed

\

In the  next result we classify the closed unit balls of $E$, $E'$ and  bounded subset $E''$ concerning if they are (or not) almost Dunford-Pettis or almost limited. As  $E'$ has the Schur property, $B_E$ is a Dunford-Pettis set. 

\begin{prop} \label{ball1}
\begin{enumerate}
\item The closed unit ball of  $E$ is not almost limited.
\item The closed unit ball of $E^{'}$ is not almost Dunford-Pettis.
\item Every norm bounded subset in $E''$ is almost limited.
\end{enumerate}
\end{prop}

\noindent \textit{Proof: } (1)
Let  $T: c_0 \to E$ and $(e_{n,n}') \subset E'$ be the positive operator  given  in the proof of Proposition \ref{wdp1}. Since $\norma{e_{n,n}'}_{B_{E}} = \sup_{x \in B_E} |e_{n,n}'(x)| \geq e_{n,n}'(Te_n) = 1$ for all $n$,where $(e_n)_n$ is the canonical basis in $c_0$, we have that $B_E$ is not almost limited.

(2) The unit diagonal sequence $e_{n,n}'' = (0, \dots, 0, e_{n}, 0, \dots)$ for all $n$ is weakly null and disjoint in $E''$. Since $\sup_{x \in B_{E'}} |e_{n,n}'' (x) | \geq \sup_{x \in B_{E'}} |e_{n,n}''(e_{n,n}')| = 1$ for all $n$, we have that $B_{E'}$ is not almost Dunford-Pettis.
As a consequence $B_{E^{'}}$ it is not almost limited.

(3)  Consider  $A \subset E''$  a norm bounded subset, then there exists $M > 0$ such that $A \subset M \cdot B_{E''} = [-Me, Me] = \sol{ Me }$.
By Lemma 2.1 of \cite{miranda}, we have that $ M \cdot B_{E''}$ is almost limited. Consequently, $A$ is almost limited as well.
\qed

\



\noindent\emph{Contact:}

\noindent M. Lilian Louren\c co, University of S\~ao Paulo, SP, Brazil, e-mail: mllouren@ime.usp.br. 

\smallskip
\noindent Vinícius  C. C. Miranda, University of S\~ao Paulo, Brazil, e-mail: vinicius.colferai.miranda@usp.br


\begin{thebibliography}{99}


\bibitem{alip} C. Aliprantis, O. Burkinshaw, \emph{Positive Operators}, Springer, Dordrecht, (2006).

\bibitem{ardakani} H. Ardakani, S. M. S. M. Mosadegh, M. Moshtaghioun and  M. Salimi, \emph{The strong Gelfand-Phillips property in Banach lattices}, Banach J. Math. Anal. 10(1) (2016) 15-26.

\bibitem{bouras} K. Bouras, \emph{ Almost Dunford-Pettis sets in Banach lattices}. Rend. Circ. Mat. Palermo 62 (2013) 227-236 .

\bibitem{bourgdiest} J. Bourgain  and J. Diestel, \emph{Limited operators and strict cosingularity}, Math. Nachrichten 119 (1984) 55-58.

\bibitem{carlo} H. Carri\'on, P. Galindo and  M. L. Louren\c co, \emph{Banach spaces whose bounded sets are bounding in bidual,} Ann. Acad. Sci. Fenn.Math. 31 (2006) 61-70.


\bibitem{car} H. Carri\'on, P. Galindo and  M. L. Louren\c co, \emph{A stronger Dunford-Pettis property}, Studia Math. 184(3) (2008) 205-216.

\bibitem{chen} J. X. Chen, Z. L. Chen, and G. X. Ji, \emph{Almost limited sets in Banach lattices}, J. Math. Anal. Appl. 412 (2014) 547-553.

\bibitem{dix} J. Diximier,\emph{ Sur un théorème de Banach, } Duke Math J. 15 (1948), 1057-1071.

\bibitem{drew1} L. Drewnowski, \emph{On Banach spaces with the Gelfand-Phillips property} Math. Z. 193 (1986) 405-411.

\bibitem{leung} D. H. Leung, \emph{On the weak Dunford-Pettis property}, Arch. Math., 52 (1989) 363-364. 

\bibitem{miranda} M. L. Louren\c co and  V. C. C. Miranda, \emph{The property (d) and the almost limited completely continuous operators,}  to apper in Anal. Math (2022), arXiv:2011.02890.




\bibitem{mach} N. Machrafi, K. El Fahri, M. Moussa and  B. Altin, \emph{A note on weak almost limited operators}, Hacet. J. Math. Stat. 48(3) (2019) 759-770.


\bibitem{meyer} P. Meyer-Nieberg, \emph{Banach Lattices}, Springer-Verlag 1991.

\bibitem{sanchez} J. A. Sanchez, \emph{Operators on Banach Lattices}, PhD thesis, Complutense University, Madrid 1985.

\bibitem{stegall} C. Stegall, \emph{Duals of Certain Spaces With the Dunford-Pettis Property}, Notices Amer. Math. Soc. 19 (1972) 799.

\bibitem{wnuk} W. Wnuk, {\em On the dual positive Schur property in Banach lattices}. Positivity 17 (2013) 759-773. 

\end{thebibliography}
\end{document}